\input amssym.def
\input amssym

\def\cz{\Bbb C}
\def\nz{\Bbb N}
\def\gz{\Bbb Z}
\def\pr{\Bbb P}

\def\r{\rightarrow}

\def\p{\partial}
\def\s{\subset}
\def\bbbc{\Bbb C}
\def\s{\subset}

\def\D{\Delta}

\def\se{\setminus}

\def\la{\lambda}
\def\ep{\epsilon}

\def\la{\lambda}

\def\v{\vartheta}

\catcode`\"=\active
\def"{\"}
\def\3{\char\ss}

\def\qed{\ifmmode\sq\else{\unskip\nobreak\hfil
\penalty50\hskip1em\null\nobreak\hfil\sq
\parfillskip=0pt\finalhyphendemerits=0\endgraf}\fi}
\def\sq{\hbox{\rlap{$\sqcap$}$\sqcup$}}
\documentstyle[12pt]{article}

\newtheorem{defi}{Definition}[section]
\newtheorem{prop}[defi]{Proposition}
\newtheorem{theo}[defi]{Theorem}
\newtheorem{lem}[defi]{Lemma}

\newtheorem{rem}[defi]{Remark}

\begin{document}

\begin{center}
{\LARGE On the Hyperbolicity of the Complements}\\
{\LARGE of Curves in  Algebraic Surfaces:}\\
{\LARGE The Three \vspace{1cm} Component Case}\\
{\large G.Dethloff, G.Schumacher, P.M.Wong}
\end{center}
\tableofcontents

\section{Introduction} In complex analysis hyperbolic manifolds have
been studied extensively, with close relationships to other areas
(cf.\ eg. \cite{LA1}). Hyperbolic manifolds are generalizations of
hyperbolic Riemann surfaces to higher dimensions. Despite the fact
that the general theory of hyperbolic manifolds is well-developed,
only very few classes of hyperbolic manifolds are known. But one could
hope that `most' of the pseudoconvex quasi-projective varieties are in
fact hyperbolic. In particular it is believed that e.g. the
complements of most hypersurfaces in $\pr_n$ are hyperbolic, if only
their degree is at least 2n+1. More precisely according to Kobayashi
\cite{KO}, and later Zaidenberg \cite{ZA} one has the following:

{\bf Conjecture:} {\it Let ${\cal C}(d_1, \ldots , d_k)$ be the space
of $k$ tupels of hypersurfaces $\, C = (C_1 , \ldots , C_k )\, $
in $\pr_n$, where ${\rm deg}(C_i)=d_i$. Then for all $(d_1, \ldots ,
d_k)$ with $\, \sum_{i=1}^k d_i =:d \geq 2n+1\, $ the set $\, {\cal
H}(d_1, \ldots , d_k)= \{ C \in {\cal C}(d_1, \ldots , d_k) : \pr_n
\se \bigcup_{i=1}^k C_i\, $ {\rm is complete hyperbolic and
hyperbolically embedded}$\}\, $ contains the complement of a proper
algebraic subset of ${\cal C}(d_1, \ldots , d_k)$.}

In this paper we shall restrict ourselves to the two dimensional case.
However we consider also more general quasi-projective complex
surfaces than the complements of curves in the projective plane.

Concerning the above conjecture, the following was known: It seems
that the conjecture is the more difficult the smaller $k$ is. Other
than in the case of 5 lines $({\cal C}(1, 1, 1, 1, 1))$, the
conjecture was previously proved by M.~Green in \cite{GRE2} in the
case of a curve $C$ consisting of one quadric and three lines
(${\cal C}(2, 1, 1, 1)$). Furthermore, it was shown for ${\cal C}(d_1,
\ldots, d_k)$, whenever $k\geq5$, by Babets in \cite{BA}. A result
which went much further was given by Eremenko and Sodin in \cite{E-S},
where they proved a Second Main Theorem of value distribution theory
in the situation $k \geq 5$. Green proved in \cite{GRE1} that for any
hypersurface $C$ consisting of at least four components in
$\pr_2$ any entire curve $f:\Bbb C \to \pr_2 \setminus C$ is
algebraically degenerate. Knowing this, it follows immediately that
for generic configurations, any such algebraically degenerate map is
constant, hence the conjecture is true for any family ${\cal C}(d_1,
\ldots, d_k)$ with $k \geq 4$ (cf. \cite{DSW}). (The degeneracy locus
of the Kobayashi pseudometric was studied by Adachi and Suzuki in
\cite{A--S1}, \cite{A--S2}).

In our paper \cite{DSW} we gave a proof of the conjecture for 3
quadrics (${\cal C}(2, 2, 2)$), based on methods from value
distribution theory. The three quadric case had been previously
studied by Grauert in \cite{GR} who used differential geometry.
However, certain technical problems still exist with this approach.
For ${\cal C}(2, 2, 1)$, i.e. two quadrics and a line, we proved with
similar methods the existence of an open set in the space of all such
configuarions, which contains a quasi-projective set of codimension
one, where the conjecture is true.

The paper contains two main results. The first is Theorem~\ref{MT}. It
states that the conjecture is true for almost all three component
cases, namely for ${\cal C}(d_1, d_2, d_3)$ with $d_1, d_2, d_3 \geq
2$ and at least one $d_i \geq 3$. Together with our result for three
quadrics (which, by the way, occur on the borderline of the method
used in this paper) this means that the conjecture is true for three
components whenever none of them is a line. We finally remark that we
get a weaker conclusion also for ${\cal C}(d_1, d_2, d_3)$ where, up
to enumeration, $d_1=1$, $d_2 \geq 3$, $d_3 \geq 4$: Namely we show that
any holomorphic map $f:\Bbb C \to X$ is algebraically degenerated,
i.e.\ $f(\Bbb C)$ is contained in a proper algebraic subset of $X$.

The other main result is Theorem~\ref{MT1}. We consider a smooth
surface $\bar X$ in $\pr_3$ of degree at least five for which every
curve on $\bar X$ is the complete intersection with another
hypersurface. Surfaces of this kind are much more general than $\pr_2$
-- by the Noether-Lefschetz theorem (cf. \cite{N-L}) the 'generic'
surface in $\pr_3$ of any given degree at least four has this property
(`generic' here indicates the complement of a countable union of
proper varieties). Let $C$ be a curve on $\bar X$ consisting of three
smooth components intersecting transversally. From our assumptions we
know that $C$ is a complete intersection of $\bar X$ and a
hypersurface $B$. We assume that the degree of $B$ is at least five.
Now Theorem~\ref{MT1} states the hyperbolicity of any such $X=\bar X
\se C$. Moreover $X$ is complete hyperbolic and hyperbolically
embedded.

Our method of proof is the following: We heavily use a theorem due to
S.~Lu (cf. \cite{Lu}). It states that for a certain class of
differentials $\sigma$, which may have logarithmic poles along the
curve $C$, and any holomorphic map $f: \cz \r \pr_2 \se C$ the
pull-back $f^*(\sigma)$ vanishes identically. This can be interpreted
as algebraic degeneracy of the tangential map corresponding to $f: \cz
\r X$. Our aim is to show algebraic degeneracy of the map $f$ itself.

The paper is organized as follows: In section~2 we collect, for the
convenience of the reader, some basics from value distribution theory.
(Readers who are familar with these may skip this section). In
section~3 we fix the notation and quote some theorems which are needed
in the following proof, especially Lu's theorem. Furthermore we
examine more closely the spaces of sections which are used in Lu's
theorem and get sections with special zero sets. The essential step of
our paper is the proof of Theorem~\ref{deg} in section~4. It states
the algebraic degeneracy of holomorphic maps $f:\Bbb C \to X$, if
${\rm Pic}(\bar X)=\Bbb Z$, and under assumptions on the determinant
bundle and the Chern numbers of the logarithmic cotangent bundle on
$\bar X$ with respect to $C$. The proof uses value distribution theory
and the existence of the special sections which were constructed in
section~3. In section~5 we compute the Chern numbers and the
determinant bundle in the situation where $\bar X$ is a complete
intersection (Theorem~\ref{main}). We apply this to $\bar X = \Bbb
P_2$, and to hypersurfaces in $\Bbb P_3$ using the Noether Lefschetz
theorem (Theorem~\ref{main2}). Finally in section~6 we apply
Theorem~\ref{main2} and get Theorem~\ref{MT}, using an argument like
in our paper \cite{DSW} to prove the nonexistence of algebraic entire
curves in generic complements. Furthermore we apply
Theorem~\ref{main2} using results of Xu \cite{Xu} and Clemens
\cite{Cl} to get Theorem~\ref{MT1}.

The first named author would like to thank S.~Kosarew (Grenoble) for
valuable discussions. The second named author would like to thank the
SFB 170 at G"ottingen, and the third named author would like to thank
the SFB 170 and the NSF for partial support.

\section{Some tools from Value Distribution Theory}

In this section we fix some notations and quote some facts from Value
Distribution Theory.  We give references but do not trace these facts
back to the original papers.

We define the characteristic function and the counting function, and
give some formulas for these.

Let $\,||z||^2= \sum_{j=0}^n |z_j|^2$, where $(z_0,\ldots ,z_n) \in
\bbbc^{n+1}$, let $\D_t = \{\xi \in \bbbc :  |\xi| < t \}$, and let $d^c =
(i/4 \pi) (\overline{\p} - \p)$.  Let $r_0$ be a fixed positive
number and let $\,r \geq r_0$.  Let $\,f:\bbbc \r \pr_n\,$ be entire, i.e.
$f$ can be written as $\, f=[f_0:\ldots :f_n]\,$ with holomorphic
functions $\, f_j :  \bbbc \r \bbbc\, , j=0,\ldots ,n\,$ without common
zeroes.  Then the {\it characteristic function} $T(f,r)$ is defined as
$$ T(f,r) = \int_{r_0}^r \frac{dt}{t} \int_{\D_t} dd^c \log ||f||^2$$
Let furthermore $\, D=V(P)$ be a divisor in $\pr_n$, given by a
homogeneous polynomial $P$.  Assume $\, f(\bbbc) \not\s \hbox{ {\rm
support}}(D)$.  Let $\,n_f(D,t)\,$ denote the number of zeroes of $\, P
\circ f\,$ inside $\, \D_t\,$ (counted with multiplicities).  Then we
define the {\it counting function} as
$$
N_f(D,r) = \int_{r_0}^r n_f(D,t) \frac{dt}{t}
$$

Stokes Theorem and transformation to polar coordinates imply (cf. \cite{WO}):
\begin{equation} \label{1}
T(f,r) =
\frac{1}{4 \pi} \int_0^{2 \pi} \log  ||f||^2 (re^{i \v})d \v + O(1).
\end{equation}

The characteristic function as defined by Nevanlinna for a holomorphic
function $\,f:  \bbbc \r \bbbc$ is
$$
T_0(f,r) = \frac{1}{2 \pi} \int_0^{2 \pi}
\log ^+ |f(re^{i \v})| d \v .
$$

For the associated map  $\, [f:1]:  \bbbc \r \pr_1$ one has
\begin{equation} \label{2}
T_0(f,r) = T([1:f],r) + O(1)
\end{equation}
(cf. \cite{HA}).

By abuse of notation we will, from now on, for a function $\, f:  \bbbc \r
\bbbc$, write $T(f,r)$ instead of $T_0(f,r)$.  Furthermore we
sometimes use $N(f,r)$ instead of $N_f({z_0=0},r)$.

We state some elementary properties of the characteristic function:

\begin{lem}\label{calc}
Let $f,g,f_j: \cz \to \cz$ be entire holomorphic functions for
$j=0,\ldots,n$. Then
\begin{description}
\item[a)]
$$
T(f\cdot g, r) \leq T(f,r) + T(g,r) + O(1)
$$
\item[b)]
$$
 T([f_0:\ldots:f_n],r)  \leq \sum_{j=0}^n T(f_j,r)
+ O(1)
$$
\item[c)]
$$
T(f+g,r) \leq T(f,r)+ T(g,r) + O(1)
$$
\end{description}
\end{lem}

{\it Proof:}\/ Propterty a) is obvious for $T_0$ and generalizes to $T$
because of (\ref{2}).

Property b) is a consequence of
$$
\log \sum_{j=0}^n|f_j|^2 \leq \sum_{j=0}^n \log (1+|f_j|^2).
$$
Property c) is a consequence of
$$
\log^+|f+g| \leq \log (1+|f|+|g|) \leq \log(1+|f|) +\log(1+|g|)
\leq \log^+|f| + \log^+|g| +2
$$
\qed

Later on we will use the concept of finite order.
\begin{defi}
Let $s(r)$ be a positive, monotonically increasing function
defined for $\,r \geq r_0$. If
$$ \overline{\lim_{r \r \infty} } \frac{\log  s(r)}{\log  r} = \la$$
then $s(r)$ is said to be of order $\la$. For entire $\,f:\bbbc \r \pr_n\,$
or $\, f: \bbbc \r \bbbc\,$ we say that $f$ is of order $\la$, if
$T(f,r)$ is.
\end{defi}

\begin{rem}\label{remfo}
Let $f=[f_0:\ldots:f_n]:\bbbc \to \pr_n$ be a
holomorphic map of finite order $\la$.  Then $\log T(f,r)= O(\log r)$.
\end{rem}

For holomorphic maps to $\pr^1$ whose characteristic function only
grows like $\log r$ we have the following characterization (cf.
\cite{HA}):
\begin{lem}\label{logr}
Let $f=[f_0:f_1]: \cz \r \pr^1$ be entire.
Then $T(f,r) =O(\log r)$ if and only if the meromorphic function
$f_0/f_1$ is equal to a quotient of two polynomials.
\end{lem}

We need the following:
\begin{lem} \label{e}
Assume that $\,f: \bbbc \r \pr_n\,$ is an entire map and misses the
divisors
$\,\{ z_j = 0\}\,$ for $j=0,\ldots,n$ (i.e. the coordinate hyperplanes
of $\pr_n$).
Assume that $f$ has order at most $\la$. Then $f$ can be written as
$\,f = [1:f_1:\ldots :f_n]\,$ with $\, f_j(\xi) = e^{P_j(\xi)}$, where
the $P_j(\xi)$
are polynomials in $\xi$ of degree $d_j\leq \lambda$.
\end{lem}
{\it Proof:} We write  $\, f=[1:f_1:\ldots :f_n]\,$ with holomorphic
$\,f_j: \bbbc \r \bbbc \se \{0\}$. Now we get with equations (\ref{1})
and (\ref{2}) for $j=1,\ldots ,n$:
$$
T(f_j,r) = T([1:f_j],r) + O(1) \leq T(f,r) + O(1),
$$
hence the $f_j$ are nonvanishing holomorphic functions of order
at most $\la$. This means that
$$
{\rm lim sup}_{r \r \infty} \frac{T(f_j,r)}{r^{\la + \ep}} =0
$$
for any $\, \ep > 0$. From this equation our assertion follows with the
Weierstra\3 theorem as it is stated in \cite{HA}. \qed

We state the First and the Second Main Theorem of Value Distribution
Theory which relate the characteristic function and the counting
function (cf.  \cite{SH}):

Let $\,f:\bbbc \r \pr_n\,$ be entire, and let $D$ be a divisor in $\pr_n$
of degree $d$, such that $\,f(\bbbc) \not\s \hbox{ {\rm support}}(D)$.
Then:

\medskip
{\bf First Main Theorem} $$ N_f(D,r) \leq d \cdot T(f,r) + O(1)$$

Assume now that $\, f(\bbbc)\,$ is not contained in any hyperplane in
$\pr_n$, and let $\, H_1,\ldots ,H_q\,$ be distinct hyperplanes in
general position.  Then

\medskip
{\bf Second Main Theorem}
$$
(q-n-1)T(f,r) \leq \sum_{j=0}^q N_f(H_j,r) + S(r)
$$
where $\: S(r) \leq O(\log (rT(f,r)))\,$ for all $\,r \geq r_0\,$ except
for a set of finite Lebesque measure.
If $f$ is of finite order, then
$\, S(r) \leq O(\log r)\,$ for all $\,r \geq r_0$.

\section{Setup and Basic Methods}

We denote by $\bar X$ a non-singular projective surface and by $C$ a
curve in $\bar X$ whose irreducible components are smooth and intersect
each other only in normal crossings. Let
$X=\bar X \setminus C$.

We denote by $E$ the dual of the bundle $\Omega^1_{\bar X}(\log C)$ of
holomorphic one forms of $\bar X$ with logarithmic poles along $C$.
Then we define the projectivized logarithmic tangent bundle $p:\pr
(E)\to \bar X$ over $\bar X$ to be the projectivized bundle whose
fibers correspond to the one dimensional subspaces of the fibers of
$E$. Furthermore let ${\cal O}_{\pr (E)}(-1)$ be the sheaf associated
to the tautological line bundle on $\pr (E)$, for which we have the
canonical isomorphism between the total space of ${\cal O}_{\pr
(E)}(-1) \se \{{\rm zero section}\}$ and the total space of $E \se \{
{\rm zero section}\}$.

Let $D$ be a divisor on $\bar X$.
According to  a Theorem of Kobayashi-Ochiai (cf. \cite{K-O}), the
cohomology, in particular the holomorphic sections, of a symmetric
power of $E^*$ tensorized with the bundle $[-D]$, corresponds to the
cohomology of the m-th power of the dual of the
tautological line bundle on $\pr(E)$,
 tensorized with the pull-back of $[-D]$:
$$
H^0(\bar X, S^m(E^*)\otimes [-D])
\simeq H^0(\pr(E),{\cal O}_{\pr(E)}(m)\otimes p^*[-D]),
$$
and in particular
$$
H^0(\bar X, S^m(E^*)) \simeq H^0(\pr(E),{\cal O}_{\pr(E)}(m)).
$$

Let $f:\cz \to X$ be a holomorphic map. Denote by
$$(f, f'): T(\cz) \r T(X)$$
the
induced map from $T(\cz)$ to the holomorphic
tangent bundle $T(X)$, which gives rise to a meromorphic map
$$F: \cz \r \pr (T(X))$$
from $\cz$
to $\pr(T(X))$. Since the domain is of dimension one, points of
indeterminacy can be eliminated, more precisely the map $F$ extends
holomorphically into the points $\xi \in \cz$ where $f' (\xi)=0$.
We denote the extended holomorphic map on $\cz$ again by $F$.
Since the restriction of
$E$ to $X$ is isomorphic to the holomorphic tangent bundle of $X$,
any  map
$$f: \cz \r X$$
has
a unique holomorphic lift $F:\cz \to \pr(E)$.

The following theorem, which is a special case of Theorem~2 of Lu in
\cite{Lu} (it actually follows already from Proposition~4.1 there)
imposes restrictions to such lifts
$F$.

\begin{theo}[Lu]\label{Lutheo}
Assume that the divisor $D$ is ample and that
there exist a non-trivial holomorphic section
$$
0\neq \sigma \in H^0(\pr(E),{\cal O}_{\pr(E)}(m)\otimes p^*[-D]).
$$
Then for any non-constant holomorphic map $f:\cz \to X$ the holomorphic
lift $F:\cz \to \pr(E)$ has values in
the zero-set of $\sigma$.
\end{theo}

In order to apply  Lu's theorem in a given situation it is
important to guarantee the existence of suitable sections.

Let $(\bar X, C)$ be given as above. The logarithmic Chern classes
$\bar c_j(X)$
are by definition the Chern classes of the logarithmic tangent bundle $E$:
$$
\bar c_j(X)= c_j(E) = c_j(X,(\Omega^1_{\bar X}(\log C))^*).
$$
Now the existence of suitable sections is guaranteed by
the following theorem of Bogomolov (cf. \cite{Bo}
and also Lu \cite{Lu} ( Proof of Proposition~3.1 and localization to
the divisor $D$).

\begin{theo}[Bogomolov]\label{Bogo} Let $D$ be a divisor on $\bar X$,
$D$ effective (i.e.\ \\ $D \geq 0$). Assume that
$$
\bar c_1^2(X) - \bar c_2(X) > 0,
$$
and that
$$
{\rm det}(E^*)
$$
is effective.
Then there exist positive
constants $A, B$ and $m_0, n_0 \in \nz$, such that
$$
A\cdot m^3 \leq h^o(\bar X, S^{mn_0}(E^*)\otimes [-D]) \leq B\cdot m^3
$$
for all $m\geq m_0$.
\end{theo}

We have the following nonexistence statement, which is a consequence of
the logarithmic version of the Bogomolov's lemma due to Sakai
(cf. \cite{Sak}). It will also become important for the following proofs.

\begin{lem}\label{lem2}
Assume that the divisor $D$ is ample. Then the following group vanishes:
$$
H^0(\pr(E), {\cal O}_{\pr(E)}(1)\otimes p^*[-D])=\{0\}.
$$
In particular, there is no logarithmic $1$-form on $\bar X$ which
vanishes on $D$.
\end{lem}

{\it Proof:}\/ The existence of a non-trivial section $s \in H^0(
\bar X, \Omega^1_{\bar X}(\log C)\otimes [-D])$ implies that the
invertible sheaf ${\cal L}:=[D]$ can be realized as a subsheaf
of $\Omega^1_{\bar X}(\log C)$. According to a result of Sakai
\cite{Sak}, (7.5), this implies that the $\cal L$-dimension of $
\bar X$ equals one, which is clearly impossible since $[D]$ is
 ample. \qed

Next we deal with divisors in $\pr (E)$ which project down to all of $\bar X$:

\begin{defi}
Consider the projection $p:\pr(E) \to \bar X$. We call a divisor
$Z\subset \pr(E)$ horizontal, if $p(Z)=\bar X$.
\end{defi}

Those horizontal divisors which occur as parts of the zero sets
$V(\sigma)$ of sections $0\neq \sigma \in H^0(\pr(E), {\cal
O}_{\pr(E)}(m)\otimes p^*[-D])$ will play an important role in the
sequel. We study this relationship somewhat closer.

\begin{lem}\label{horcomp}
Given
$$
0\neq \sigma \in H^0(\pr (E), {\cal O}(m)\otimes p^*[-D])
$$
there exist divisors $E_j$ $ j=1, \ldots, l$ on $\bar X$, and numbers
$a_j, n_j\in \nz$ such that $[\sum a_j\cdot E_j - D] \geq
0$ and sections $s_j\in H^0(\pr(E), {\cal O}(n_j)\otimes p^*[-E_j])$,
$\tau \in H^0(\pr(E), p^*[\sum a_j\cdot E_j - D])$ such that $\sigma =
\tau \otimes_{1\leq j \leq l} s_j^{a_j}$ with the following property:
The zero-sets of $s_j$ are precisely the irreducible horizontal
components of $V(\sigma)$. \end{lem}

{\it Proof:}\/ Let $0\neq \sigma \in H^0(\pr(E), {\cal O}(m)\otimes
p^*[-D])$ be a non-trivial section and $V(\sigma)$ its zero divisor.
We denote by $S_j$; $j=1, \ldots l$ the irreducible horizontal
components of $V(\sigma)$. Since ${\rm Pic }(\pr(E))={\rm Pic }(X)
\oplus \gz$, we get $[S_j]= {\cal O}_{\pr(E)}(n_j)\otimes p^*[-E_j]$
for certain divisors $E_j\subset \bar X$ with $n_j \geq 1$. This fact
follows by restricting the bundles $[S_j]$ to a generic fiber of $p$.
Let $a_j$ be the multiplicities of $\sigma$ with respect to $S_j$,
then in particular $a_1n_1+\ldots a_ln_l =m$. (This fact follows again
by restricting bundles and sections to a generic fiber of $p$.)
Canonical sections of $[S_j]$ give rise to non-trivial sections
$s_j\in H^0(\pr(E), {\cal O}(n_j)\otimes p^*[-E_j])$ which vanish
exactly on $S_j$. Thus $\tau:= \sigma/( s_1^{a_1}\cdot\ldots\cdot
s_l^{a_l})$ is a (holomorphic) section of $H^0(\pr(E), p^*[\sum
a_j\cdot E_j - D])$. In particular $[\sum a_j\cdot E_j - D] \geq 0$.
\qed

In order to control the horizontal divisors of a section of
$H^0(\pr(E), {\cal O}(m)\otimes p^*[-D])$, the number $m$ will be
chosen minimal in the following sense.

For any $k\in \nz$ we set
$$
\mu_k:=\inf\{m;h^0(\pr(E), {\cal O}(m)\otimes p^*[-kD])\neq 0 \}.
$$
and
$$
\mu := \inf_{k \in \nz} \{\mu_k\}.
$$

\begin{lem}\label{onecomp}
Assume that ${\rm Pic}(\bar X)= \gz$ and that $[D]$ is the ample generator
of ${\rm Pic}( \bar X)$.\\
Then we have
$$
2 \leq \mu < \infty
$$
and if $k_0 ={\rm min}\{k \in \nz :\mu_k = \mu\}$, there exists a
non-trival section
$$
0\neq \sigma \in H^0(\pr (E), {\cal O}(\mu)\otimes p^*[-k_0D])
$$
such that exactly one horizontal
component of $V(\sigma)$ exists and has multiplicity one.
\end{lem}

{\it Proof:}\/ Since some multiple of $D$ is a very ample and hence
linear equivalent to an effective divisor, we have $\mu < \infty$ from
Theorem~\ref{Bogo}. From Lemma~\ref{lem2} we then get $\mu \geq 2$.\\
Now take any section $0\neq \sigma \in H^0(\pr (E), {\cal
O}(\mu)\otimes p^*[-k_0D])$. We use Lemma~\ref{horcomp}. Since $[D]$
is a generator of ${\rm Pic }(\bar X)$, there exist $b_j \in \gz$ such
that $[E_j]=b_j\cdot [D]$. Since $[\sum a_j\cdot E_j - k_0 D] \geq 0$
we have $\sum a_jb_j \geq k_0$. Since all $a_j\geq 0$, there must be
at least one $b_j>0$, say $b_1> 0$. Now $s_1\in H^0(\pr(E), {\cal
O}(n_1)\otimes p^*[-b_1\cdot D])$ is a non-trivial section. By
definition of $\mu$ we have $n_1\geq \mu$ which means $n_1=\mu$, since
$\sum a_jn_j = \mu$. So in terms of the notion of Lemma~\ref{horcomp}
$\sigma=\tau \cdot s_1$, i.e.\ $S$ contains only one horizontal
component. This component has multiplicity one. \qed

\section{Algebraic Degeneracy Of Entire Curves}

Let  $\bar X $ be a non-singular (connected) projective surface.

\begin{defi} Let $f:\cz \to \bar X$ be a holomorphic map. We call
$f$ {\em algebraically degenerate}, if there exists an algebraic curve
$A \subset \bar X$ such that $f(\cz)$ is contained in $A$.
\end{defi}

Our main result on algebraic degeneration is:

\begin{theo}\label{deg}
Let $C \subset \bar X$ be a curve consisting of three
smooth components with normal crossings.
Assume that:\\
i) ${\rm Pic}(\bar X) = \gz$\\
ii) The logarithmic Chern numbers of $X = \bar X \se C$
satisfy the inequality $$\bar c_1^2(X) - \bar c_2(X) > 0$$
iii) The line bundle ${\rm det}(E^*)$ is effective, where
$E^* = \Omega_X^1(\log  C)$ is the logarithmic cotangent bundle.\\
Then any holomorphic map $f:\cz \to \bar X\setminus C$ of order at most two is
algebraically degenerate.
\end{theo}

Remark: The theorem also holds without the assumption on the order of
the map $f$, but since we are mostly interested in the hyperbolicity
of the complement, we include this assumption, because it slightly
simplifies the proof.\\

The rest of this section is devoted to the proof of this Theorem.

Let again $[D]$ be an ample generator of ${\rm Pic}(\bar X)$. Let $k
\in \nz$ be a natural number such that $[kD]$ is very ample. Then by
Theorem~\ref{Bogo} there exists a symmetric differential $\omega \in
H^0(\bar X, S^m(E^*) \otimes [-kD])$ which is not identically zero.
By Theorem~\ref{Lutheo} we know that $f^*\omega \equiv 0$.\\

The proof now will work as follows: The three components of the curve
$C$ give rise to a morphism $\Phi:\bar X \to \pr_2$ which maps $C$ to
the union of the three coordinate axis. In the first step of the proof
we show that we can `push down' the symmetric differential $\omega$
by this morphism to some symmetric rational differential $\Omega$ on
$\pr_2$ and that we still have $(\Phi \circ f)^*(\Omega) \equiv 0$.
Since $\Phi \circ f$ maps the complex plane to the complement of the
three coordinate hyperplanes in $\pr_2$, we will be able to interpret
this, in the second step of the proof, as an equation for nonvanishing
functions with coefficients which may have zeroes, but which grow of
smaller order, only. In such a situation we then can apply Value
Distribution Theory.


{\bf First step:} We first remark that the intersection number of any
two curves $D_1$ and $D_2$ is positive (including self intersection
numbers). Let $[D]$ be the ample generator of ${\rm Pic}(\bar X) =
\gz$. Now $[D_j] = a_j [D]$; $a_j \in \gz$, and $0< D_j \cdot D = a_j
D^2$ (cf.\ the easy implication of the Nakai criterion). Hence all
$a_j$ are positive, and
$$
D_1 \cdot D_2 = a_1 a_2 D^2 >0.
$$

We can find $a_j \in \Bbb N$; $ j=1, 2, 3$ such that
$[a_1C_1]=[a_2C_2]=[a_3C_3]$, since the divisors $C_j$ $ j=1, 2, 3$
are effective. Let $\sigma_j \in H^0(\bar X, L)$ be holomorphic
sections which vanish exactly on $C_j$. Then
$$
\Phi = [\sigma_1:\sigma_2:\sigma_3]: \bar X \r \pr_2
$$
defines a rational map, which is a morphism, since the three components
do not pass through any point of $\bar X$.

\begin{lem} \label{dom}
The morphism $\Phi$ is a branched covering.
\end{lem}

{\it Proof:}\/ Since $C_2\cdot C_3 >0$, the fiber $\Phi^{-1}(1:0:0) =
C_2 \cap C_3$ is non-empty. By assumption $C_2 \cap C_3$ consists of
at most finitely many points. Hence $\Phi$ is surjective and has
discrete generic fibers. Finally $\Phi$ has no positive dimensional
fibers at all: Applying Stein factorization we would get a
bimeromorphic map. Since there are no curves of negative
self-intersection, no exceptional curves exist on $\bar X$(cf.
\cite{BPV}). Hence there exist no positive dimensional fibers of
$\Phi$. \qed

Hence the morphism $\Phi$ is a finite branched covering of
$\bar X$ over $\pr_2$ with, let us say $N$ sheets. Let $R$ be
the ramification divisor of $\Phi$, $B=\Phi(R)$ the branching locus
and $R'=\Phi^{-1}(B)$. Then
$$\Phi : \bar X \se R' \r \pr_2 \se B$$
is an unbranched covering with $N$ sheets.

We now want to construct a meromorphic symmetric $mN$-form $\Omega$
defined on $\pr_2 \se B$ from the meromorphic symmetric $m$-form
$\omega$ on $\bar X$: For any point $w^0 \in \Bbb P_2 \se B$, there
exists a neighborhood $U=U(w^0)$ of $w^0$ and $N$ holomorphic maps
$a_i(w)$, $a_i: U\to \bar X \se R'$; $ i=1, \ldots, N$ such that $\Phi
\circ a_i= {\rm id}_U$. By pulling back the symmetric $m$-form
$\omega$ by means of these maps we get $N$ meromorphic symmetric
$m$-forms $ (a_i)^*(\omega)(w)$ on $U$. Taking now the symmetric
product of these $m$-forms, we get the symmetric $(Nm)$-form $\Omega$
on $U$:

$$\Omega (w) = \prod_{i=1}^N a_i^*\omega(w). $$
Let $M=Nm$.
Defining $g=\Phi \circ f$, we then have:

\begin{lem} \label{meromext}
The form
$\Omega$ extends to a rational symmetric $M$-form on $\pr_2$, which we
again denote by $\Omega$. We have $\Omega \not\equiv 0$, but
$g^*\Omega =0$.
\end{lem}

The first statement of this lemma is probably well known, and the
second statement is considered to be obvious. But since we did not
find a reference, we will include a proof of this Lemma at the end of
this section.

We proceed with the proof of Theorem~\ref{deg}.

Denote the homogeneous coordinates of $\pr_2$ by $w_0,w_1,w_2$. On
$\pr_2 \se V(w_0)$ we have inhomogeneous coordinates
$\xi_1=w_1/w_0, \xi_2 = w_2/w_0$. Hence on $\pr_2 \se
V(w_0)$ the symmetric $M$-form $\Omega$ can be written as
\begin{equation}\label{om}
\Omega = \sum_{i=1}^M R_i(\xi_1,\xi_2) (d\xi_1)^i (d\xi_2)^{M-i}
\end{equation}
where multiplication means the symmetric tensor product here, and the
coefficients $R_i(\xi_1, \xi_2)$ are rational functions in $\xi_1$ and
$\xi_2$.

Now $g:\cz \r \pr_2$ has values in the complement $\pr_2 \se
V(w_0w_1w_2)$ of the three coordinate axis, hence the functions $g_j
=\xi_j \circ g$ are holomorphic and without zeroes. Since $g^*\omega
\equiv 0$ on $\cz$, equation (\ref{om}) implies
\begin{equation}\label{omg}
\sum_{i=1}^M R_i(g_1(\eta),g_2(\eta)) (g_1'(\eta))^i(g_2'(\eta))^{M-i}
\equiv 0
\end{equation}
for all $\eta \in \cz$. This equation still holds if we clear the
denominators of the $R_i(\xi_1, \xi_2)$ simultaneously, so without
loss of generality we may assume from now on that in equation
(\ref{omg}) the $R_i(g_1(\eta), g_2(\eta))$ are polynomials in
$g_1(\eta)$ and $g_2(\eta)$, i.e. we have
\begin{equation} \label{terms}
R_i(g_1(\eta),g_2(\eta))= \sum_{j,k} a_{ijk} (g_1(\eta))^j (g_2(\eta))^k
\end{equation}

Under our assumptions we are able to say more about the functions
$g_i$; $i=1, 2$. Since the holomorphic map $f:\cz \r \bar X \se C$ was
of finite order at most two, this is also true for $g=\Phi \circ f$ by
Lemma~\ref{calc}, since the components of $g$ are polynomials in the
components of $f$.

Hence by Lemma~\ref{e}, we have
\begin{equation} \label{gi}
g_i (\eta) = \exp (p_i(\eta))
\end{equation}
where the $p_i(\eta)$; $ i=1,2$ are polynomials in $\eta$ of degree at most
two.
Furthermore we may assume that both polynomials are non-constant, otherwise
$g$ would be linearly degenerate and so $f$ would be algebraically
degenerate, and we were done.

Replacing equation (\ref{gi}) and equation (\ref{terms})
in equation (\ref{omg}) we get
\begin{equation}\label{new}
\sum_{i=1}^M \sum_{j,k} a_{ijk} \exp\{(i+j)p_1(\eta))+(M-i+k)p_2(\eta)\}
   (p_1'(\eta))^i(p_2'(\eta))^{M-i}
\equiv 0.
\end{equation}
If we still allow linear combinations of the above summands with
constant coefficients $c_{ijk}$ in
equation (\ref{new}) we can pass to a subset $S$ of indices which occur in this
equation and get a relation
\begin{equation}\label{news}
\sum_{(i, j, k) \in S} c_{ijk} a_{ijk}
\exp\{(i+j)p_1(\eta))+(M-i+k)p_2(\eta)\}
(p_1'(\eta))^i(p_2'(\eta))^{M-i} \equiv 0
\end{equation}
but now with the additional property that $S$ is minimal with
equation (\ref{news}).
Let $S$ have $L$ elements.

Since we may assume that the polynomials $p_i(\eta)$ are nonconstant,
we know that the $p_i' (\eta)$ are not identically zero and hence that
$L \geq 2$.

For the rest of this proof
we will distinguish between two cases:\\
{\bf Case 1:}\/ There exist two summands in equation (\ref{news}) the quotient
 of which is not a rational function in the variable $\eta$.\\
{\bf Case 2:}\/
The quotient of any two summands in equation (\ref{news}) is a rational
function in the variable $\eta$.\\
We shall show that the first case is impossible whereas in the second
case algebraic degeneracy is shown.

{\bf Case 1:} We could immediately finish up the proof under the
assumptions of case 1 by using a Second Main Theorem for moving
targets to equation~(\ref{news}), as to be found e.g. in the paper of
Ru and Stoll \cite{R-S}. Another approach is to treat
equation~(\ref{news}) directly with a generalized Borel's theorem (we
can regard this equation as a sum of nonvanishing holomorphic
functions with coefficients which may vanish, but which grow of a
smaller order than the nonvanishing functions, only). We present here
a more elementary argument based on the Second Main Theorem which
might also be considered somewhat simpler.

First, it is easy to see that $L \geq 3$, since for $L=2$ we would get,
by dividing in equation~(\ref{news}) through one of the exponential
terms, that the exponential of a nonconstant polynomial is equal to a
quotient of two other polynomials, which is absurd.

Let $\psi_1,\ldots, \psi_L$ be some enumeration of the summands which
occur in equation (\ref{news}). Then, after factoring out possible
common zeroes of the entire holomorphic functions $\psi_1, \ldots,
\psi_L$ we get an entire holomorphic curve
$$
\Psi : \cz \r \pr^{L-1};
\eta \r [\psi_1(\eta):\ldots:\psi_L(\eta].
$$

If we denote the homogenous coordinates of this $\pr^{L-1}$ by
$[z_1:\ldots:z_L]$, the image of $\Psi$ is contained in the hyperplane
$H=\{z_1+\ldots+z_L=0\}$ and does not hit any of the coordinate
hyperplanes $H_i=\{z_i=0\}$. So we can regard $\Psi$ also as an entire
holomorphic mapping with values in the hyperplane $H$ (which is
isomorphic to $\pr^{L-2}$) which does not intersect the $L$ different
hyperplanes $H \cap H_i$ in $H$. It is now an important fact that
these hyperplanes are in general position in $H$, and that the entire
curve $\Psi$ is not mapping $\cz$ entirely into any hyperplane in $H$
(the latter follows from the minimality condition in equation
(\ref{news})), because under these conditions we can apply the Second
Main Theorem (cf. section~2), which yields:
\begin{equation}
\label{smt} (L-(L-2)-1) T(\Psi, r) \leq \sum_{i=1}^L N_{\Psi}(H \cap
H_i, r) + O(\log r)
\end{equation}
because the entire curve $\Psi$ is
of finite order at most two by Lemma~\ref{calc}. Now we have
$$
N_{\Psi}(H \cap H_i, r) = N(\{\psi_i=0\}, r) \leq M (N(\{p_1' =0\},
r)+N(\{p_2' =0\}, r)).
$$
The First Main Theorem (cf. section~2) and
Lemma~\ref{logr} imply that the right hand side grows at most of order
$O(\log r)$ only, so equation (\ref{smt}) yields that
\begin{equation} \label{log}
T(\Psi, r) =O(\log r).
\end{equation}
We know by the
assumption of case 1 that there exist indices $i, j$ such that
$\psi_i(\eta)/\psi_j(\eta)$ is not a rational function in $\eta$. Then
$[\psi_i(\eta):\psi_j(\eta)]:\cz \r \pr^1$ is an entire curve for
which by Lemma~\ref{logr} the characteristic function
$T([\psi_i:\psi_j])$ grows faster than $\log r$, so (by the formula
for the characteristic function given in equation~(\ref{1})) this is
also true for $T(\Psi, r)$ contradicting equation~(\ref{log}). So we
have shown that under the assumptions of case 1 we get a
contradiction.

{\bf Case 2:} We want to show first that there exist nonvanishing
complex numbers $\gamma$ and $\lambda$ such that
\begin{equation} \label{deriv}
\lambda p_1'(\eta) = \gamma p_2'(\eta).
\end{equation}

We only need to show that $p_1'(\eta)$ and $p_2'(\eta)$ are linearly
dependent, because if one of them is the zero polynomial, we have
algebraic degeneracy of $g$ and hence of $f$. So assume that
$p_1'(\eta)$ and $p_2'(\eta)$ are linearly independent. Then no linear
combination of $p_1(\eta)$ and $p_2(\eta)$ is a constant polynomial.
So under the assumptions of case 2 get that for all $(i, j, k) \in S$
the terms $i+j$ in the summands
$$ c_{ijk} a_{ijk} \exp((i+j)p_1(\eta))+(M-i+k)p_2(\eta))
   (p_1'(\eta))^i(p_2'(\eta))^{M-i}
$$
are equal, and also the terms $k+(M-i)$ are the same as well. But then
for a given $i_0$ there can be at most one $(i_0, j, k) \in S$. So by
factoring out the exponential function in equation (\ref{news}) we get
a nontrivial homogenous equation of degree $M$ in $p_1'(\eta)$ and
$p_2'(\eta)$, which then can be factored in linear factors. Since then
one of the linear factors has to vanish identically we get the linear
dependency of $p_1'(\eta)$ and $p_2' (\eta)$ again, so the assumption
of linear independency was wrong.

We now want to construct a special symmetric form with at most
logarithmic poles as singularities along the curve $C$ which is
annihilated by $f$.

Let us simply state equation~(\ref{deriv})
in terms of the original entire curve $f$.  We have
\begin{equation} \label{ww}
p_i'(\eta) = \frac{dg_i(\eta)}{g_i(\eta)} =
 (\Phi \circ f)^* \frac{d\xi_i}{\xi_i}=
f^* \omega_i
 \end{equation}
where $\omega_i$; $i=1,2$ is a differential one form on $\bar X$
with at most logarithmic
poles along $C$. Define $\omega_0 = \lambda \omega_1 - \gamma \omega_2$.
Then $\omega_0 \in H^0(\bar X, E^*)$, and since
$$\omega_0 =\Phi^* (\lambda \frac{d\xi_1}{\xi_1} - \gamma
\frac{d\xi_2}{\xi_2}) $$
and the map $\Phi$ is a local isomorphism outside the branching, we have
$$\omega_0 \not= 0$$
Furthermore by equations (\ref{ww}) and (\ref{deriv}) we have
$$ f^*\omega_0  \equiv 0$$

Now the proof of the fact that $f$ is algebraically degenerate is almost
finished:

Let $\sigma \in H^0(\pr(E), {\cal O}(\mu) \otimes p^*[-k_0D])$ be the section
constructed in Lem\-ma~\ref{onecomp} and $\tilde\sigma \in H^0(\pr(E), {\cal O}
(1))$ the section which corresponds to $\omega_0$.
We recall that both sections are nontrivial, that $\mu \geq 2$, and that
$V(\sigma)$
contains only one horizontal component, which we will denote by $S_{\sigma}$,
with multiplicity one. We also recall that by Theorem~\ref{Lutheo}, the lift of
$f$ to $\pr(E)$, which we denoted by $F$,
 maps entirely into $V(\sigma)$. We may
assume that it maps into $S_{\sigma}$, otherwise by projecting down to
$\bar X$ we get that $f$ is algebraically degenerate and we are done.

If $\tilde\sigma$ does not vanish identically on $S_{\sigma}$, $F$ maps
into the zero set of $\tilde\sigma$ in $S_{\sigma}$, which has codimension
at least two. So projecting down to $\bar X$ again yields algebraic
degeneracy of $f$.

Hence we now may assume that $\tilde\sigma$ vanishes identically
on $S_{\sigma}$. Since $\tilde\sigma \in H^0(\pr(E), {\cal O}(1))$ the
degree of $V(\tilde\sigma)$ with respect to a generic fiber of the map
$p:\pr(E) \r \bar X$ is one (cf. the argument in the proof of
Lemma~\ref{horcomp}). However since $S_{\sigma}$ is the only
horizontal component of the zero set of $\sigma \in H^0(\pr(E),
{\cal O}(\mu) \otimes p^*[-k_0D])$ with $\mu \geq 2$ and has multiplicity one,
and since $\tilde\sigma$ vanishes on $S_{\sigma}$, the degree of
$V(\tilde\sigma)$ with respect to such a generic fiber must be at least
two, which is a contradiction.
 So this case cannot occur and the proof of
Theorem~\ref{deg} is complete. \qed

\begin{small}
{\it Proof of Lemma~\ref{meromext}:}\/ The assertion $g^*\Omega \equiv
0$ is clear from $f^*\omega \equiv 0$ and the construction of
$\Omega$.

In order to prove the assertion $\Omega \not\equiv 0$,
we choose a point $\xi^0 \in \pr_2 \se (B \cup \{w_0=0\})$.
In a small neighborhood $U(\xi^0)$ we have the $N$ biholomorphic
functions $a_i(\xi), i=1,...,N$ which invert the map $\Phi$ on $U(\xi^0)$.
Then we have
\begin{equation} \label{express}
((a_i)^*(\omega))(\xi) = \sum_{j=0}^m b_{ij}(\xi)
(d\xi_1)^j(d\xi_2)^{m-j}
\end{equation}
After possibly moving the point $\xi^0$ in $U(\xi^0)$ we may assume
that the meromorphic functions $b_{ij}(\xi)$ either vanish
identically on $U(\xi^0)$ or have no zero
or singularity in $\xi^0$. Let now for each $i=1,...,N$ the index
$j(i)$ be the maximal $j \in \{0,...,m\}$ such that $b_{ij}(\xi^0) \not=
0$. Let $k=\sum_{i=1}^N j(i)$. Then the $(d\xi_1)^k(d\xi_2)^{M-k}$-monomial
of $\Omega$ in the point $\xi_0$ is equal to $\prod_{i=1}^N b_{ij(i)}(\xi^0)$,
which is not equal to zero by construction.

Last we have to show that $\Omega$ extends to a rational symmetric
$M$-form on $\pr_2$.

 We only have to show how $\Omega$ can be extended over smooth points of the
branching locus, because by Levi's extension theorem (cf. \cite{G--R}) we
then can extend it over the singular locus  which is of codimension
two. Then, by Chow's Theorem it is rational.
So assume $P \in B$ is a smooth point of $B$. Then (cf. \cite{G-R})
there exists a neighborhood of $P$ over which $\Phi$ is an analytically
branched covering of a very special form:
For every point $Q$ over $P$ one can introduce local coordinates
$\xi_1,\xi_2$ around $P$ and $z_1,z_2$ around $Q$ such that
$\xi_1(P)=\xi_2(P)=z_1(Q)=z_2(Q)=0$, and neighborhoods
$U=\{ |\xi_1|<1, |\xi_2|<1\}$, $V= \{|z_1|<1, |z_2|<1 \}$ such that,
for some $b \in \{1,...,N\}$, we have
\begin{equation} \label{abc}
\Phi :V \r U; (z_1,z_2) \r (z_1^b,z_2)
\end{equation}
In order to prove our assertion in a neighborhood of $P \in B$, it is
sufficient
to prove it for the analytically branched covering in equation (\ref{abc}).

For $k=0,...,b-1$ let
$$g_k: V \r V; (z_1,z_2) \r (\exp(\frac{2\pi i k}{b}) z_1,z_2).$$
Then $G= \{g_0,...,g_{b-1} \}$ is just the group of deck transformations,
i.e.\ automorphisms which respect $\Phi$.
For the meromorphic symmetric $m$-form $\omega$ on $V$, let $\tilde\Omega$
be the symmetric product of the $b$ meromorphic symmetric $m$-forms
$(g_i)^*(\omega)$ on $V$. We are done, if we show that by projecting
down with $\Phi$ this form gives rise to a meromorphic symmetric
$M$-form on $U$, because in $U \se B$ this is just the form $\Omega$.
The symmetric $M$-form $\tilde\Omega$ can be uniquely written in the form
\begin{equation} \label{uni}
\tilde\Omega (z_1,z_2) =
\sum_{i=0}^M r_{ij}(z_1, z_2)(\frac{dz_1}{z_1})^i
(dz_2)^{M-i}
\end{equation}
with meromorphic functions $r_{ij}$ in the variables $(z_1,z_2)$.
Now $\tilde\Omega(z_1, z_2)$ is invariant under the action of $G$,
$(\frac{dz_1}{z_1})^i$ and $(dz_2)^{M-i}$ are also $G$-invariant.
Moreover $b \frac{dz_1}{z_1}=\frac{d\xi_1}{\xi_1}$ and $dz_2=d\xi_2$.
Hence the $r_{ij}$ are $G$-invariant functions, i.e.\ these are
pull-backs of meromorphic functions on $U$. \qed
\end{small}

\section{Application to the Projective Plane and Complete Intersections}
We shall apply Theorem~\ref{deg}. Throughout this section, we make the
following assumptions:

Let the smooth complex surface $\bar X$ be a complete intersection
$$\bar X = V_2^{(a_1, \ldots , a_r)} \s \pr_{r+2}\, , \: r \geq 0 $$
of hypersurfaces of degrees $a_j$, $j=1, \ldots , r$ in $\pr_{r+2}$.
Set $A=\prod_{i=1}^r a_i$ and $a=\sum_{i=1}^r a_i$. Let smooth curves
$C_j$, $j=1, 2, 3$ be given in $\bar X$ which intersect in normal
crossings. We assume that these curves are transversal intersections
of $\bar X$ with hypersurfaces of degrees $b_j$. We set
$b=b_1+b_2+b_3$.

\begin{lem}\label{riro}
 The Euler numbers of $\bar X$ and $C_j$ are:
$$e(\bar X)=A(2+(a-r-1)^2)$$
and
$$e(C_j)=Ab_j(3+r-a-b_j).$$
\end{lem}

The {\it Proof}\/ is a direct consequence of the Riemann-Roch Theorem.
According to \cite{HI}, Theorem~22.1.1, the $\chi_y$-characteristic
of a complete intersection can be computed from a generating function.
Its value at $y=-1$ yields the Euler number. \qed

In order to determine, when the assumptions of Theorem~\ref{deg} are
satisfied, we first compute $\bar c_1^2(X) - \bar c_2(X)$.

\begin{prop}\label{num}
In the above situation
$$
\bar c_1^2(X) - \bar c_2(X) = A((a-r-3)(b-4) -6 + \sum_{i<j}b_i\cdot b_j),
$$
and
$${\rm det}(E^*)= (a+b-3-r)\tilde{H},$$
where $\tilde{H}$ is a hyperplane section.
\end{prop}

{\it Proof:}\/ According to a result of Sakai
\cite{Sak} we have ${\rm det}(E^*)=[\Gamma]$, where $\Gamma= K_{\bar
X} + C$.
Then the second claim follows from the Adjunction Formula.
 Furthermore (cf. \cite{Sak}),
$$
c_1^2(E)-c_2(E)=
c_1^2(E^*)-c_2(E^*)= \Gamma^2 - e(\bar X \setminus C) =  \Gamma^2
- e(\bar X) +  e(C),
$$
where $\Gamma^2$ denotes the self intersection. It equals
$$
\Gamma^2= A(a+b-r-3)^2.
$$
For the Euler number of $\bar X$ we use Proposition~\ref{riro}. The
Euler number of $C$ is evaluated in terms of the Euler numbers
$e(C_j)$ of the
components  and the respective intersection numbers
$$
C_i\cdot C_j= A b_i b_j
$$
to be
$$
e(C) = \sum_{j=1}^3 e(C_j) - \sum_{i<j} C_i\cdot C_j.
$$
{}From these equalities we get immediately the above formula for
$c_1^2-c_2$.\qed

Now Theorem \ref{deg} yields

\pagebreak

\begin{theo}\label{main}
Let $X=\bar X \setminus C$ as above. Then any entire holomorphic curve
$f:\cz \to X$ of order at most two is algebraically degenerate, if
\begin{itemize}
\item[i)]
${\rm Pic}(\bar X)= \gz$
\item[ii)]
$(a-r-3)(b-4)+ \sum_{i<j}b_ib_j > 6$
\item[iii)]
$a+b \geq r+3$
\end{itemize}
\end{theo}
The {\it Proof}\/ follows from Theorem~\ref{deg}, and
Proposition~\ref{num}. \qed

\begin{theo}
\label{main2}
Let $X=\bar X \setminus C$ as above. Then any entire holomorphic curve
of order at most two
$f:\cz \to X$ is algebraically degenerate in any of the following cases:
\begin{itemize}
\item[a)]
${\rm Pic}(\bar X)= \gz$,
and $a \geq r+3$, $b \geq 5$.
\item[b)]
$\bar X \s \pr_3$ is a `generic' hypersurface of degree at least four,
and $b \geq 5$.
\item[c)]
Let $\bar X =\pr_2$  (i.e.\ $a_1= \ldots a_r =1$, $r \geq 0$).
Let $b_1, b_2, b_3 \geq 2$ and at least one  $b_j \geq 3$,
or up to enumeration
$b_1=1, b_2\geq 3, b_3 \geq 4$.
\end{itemize}
\end{theo}
Remark: `Generic' indicates the complement of a countable union of
proper varieties in space of all hypersurfaces.

{\it Proof:}\/ Case a) is obvious. Case b) is an application of the
Noether-Lefschetz theorem \cite{N-L} and case a).
For case c) we set e.g.\ $r=a_1=1$. Then
$$
\bar c_1^2(X) - \bar c_2(X) = -3(b-4) -6 + \sum_{i<j}b_i\cdot b_j
$$
is equal to
$$
(b_1-2)(b_2-2)+(b_1-2)(b_3-2)+(b_2-2)(b_3-2)+b-6$$
or to
$$ (b_1-1)(b_2-1)+(b_1-1)(b_3-2)+(b_2-3)(b_3-4)+(2b_2+b_3)-9$$
 From these facts the assertion of case c) follows immediately.\qed

\section{Algebraic Degeneracy of Entire Curves Versus Hyperbolicity}

\begin{theo} \label{MT}
Let $C$ be the union of
three smooth curves $C_j$ $ j=1,2,3$ in $\pr_2$ of degree $d_j$
with
$$
d_1,d_2,d_3\geq 2 \hbox{ and at least one } d_j\geq 3.
$$
Then for generic such configurations
 $ \pr^2 \se C$ is complete hyperbolic and hyperbolically
embedded in $\pr_2$.

More precisely this is the case, if the curves intersect only in normal
crossings, and if
one curve is a quadric
there must not exist a line which intersects
the two other curves only in one point each and which intersects the
quadric just in these two points.\\
\end{theo}

{\it Proof:}\/
In order to prove that
$\pr_2 \se C$ is hyperbolic and hyperbolically embedded in $\pr_2$,
we only will have to prove,
by an easy Corollary of a Theorem of M.Green (cf. \cite{GRE2}),
that there does not exist a non-constant entire curve
$\, f:\bbbc \r \pr_2 \se C$ of order at most two.

We know from Theorem~\ref{main2} that the entire curve $f:\cz
\to \pr_2\setminus C$ of order at most two
 is contained in an algebraic curve $A\subset
\pr_2$ of degree $d_0$ say.

Now the proof is almost the same as in \cite{DSW}.
Assume that there exists an irreducible algebraic curve $A \s \pr_2 $
such that $A\setminus C$ is not hyperbolic.
We know that $  A \cap C$ consists of at least
2 points $P$ and $Q$. Moreover, $A$ cannot have a singularity at $P$
or $Q$ with different tangents, because $A$ had to be reducible in such
a point, and $A\setminus C$ could be identified with an
irreducible curve with
at least three punctures. (This follows from blowing up such a point or
considering the normalization).

So $ A \cap C$ consists of exactly 2 points $P$ and $Q$
with simple tangents.  We denote the multiplicities of $A$ in $P$ and $Q$ by $
m_P $ and
$ m_Q $.   Then the inequality (cf.  \cite{FU})
$$
m_P(m_P-1)+m_Q(m_Q-1) \leq (d_0-1)(d_0-2)
$$
implies
\begin{equation} \label{*}
m_P , m_Q < d_0 \hbox{ {\rm or} }d_0=m_P=m_Q=1.
\end{equation}

  After a suitable
enumeration of its components we may assume that $P \in C_1
\cap C_2 $ and $ Q \in C_3$. If $Q \not\in C_2
\cup C_1 $ we are done, since then we may assume that $A$ is
not tangential to $ C_2$, and then, computing intersection multiplicities
according to \cite{HA}, we have
$$
m_P = I(P,A \cap C_2) = d_2d_0
$$
which contradicts equation (\ref{*}).
So we may assume that $\, Q \in C_2 \cap C_3$. Now $A$ has to
be tangential to $C_1$ in $P$ and to $\, C_3$ in $Q$,
otherwise we again get $\, m_P=d_1d_0\,$ or $\,m_Q =d_3d_0\,$ what
contradicts equation (\ref{*}). But then $\,C_2$ is not tangential
to $A$ in $P$ or $Q$, so we have
$$m_P + m_Q = I(P,A \cap C_2) + I(Q, A \cap C_2) = d_2d_0$$
Again by equation (\ref{*}) this is only possible if $d_2=2$ and
$\, m_P=m_Q=d_0=1$, but then we are in a situation which we excluded
in Theorem~\ref{MT}, which is a contradiction. \qed

We make the same assumptions as in section~5.

\begin{theo} \label{MT1}
Let $\bar X \s \pr_3$ be a `generic' smooth hypersurface of degree $d
\geq 5$ and $b \geq 5$. Then $X=\bar X \se C$ is hyperbolic and
hyperbolically embedded in $\bar X$.
\end{theo}

{\it Proof:}\/ According to Xu \cite{Xu} and Clemens \cite{Cl} $\bar
X$ does not contain any rational or elliptic curves. Hence
Theorem~\ref{main2} yields the claim. \qed


\pagebreak

\bigskip
\bigskip

\noindent Gerd Dethloff\\
Mathematisches Institut der Universit\"at G\"ottingen\\
Bunsenstra\3e 3-5\\
3400 G\"ottingen\\
Germany\\
\vspace{0.8cm}e-mail: DETHLOFF@CFGAUSS.UNI-MATH.GWDG.DE\\

\noindent Georg Schumacher\\
Ruhr-Universit\"at Bochum, Fakult\"at f\"ur Mathematik\\
Universit\"atsstra\3e 150\\
4630 Bochum 1\\
Germany\\
\vspace{0.8cm}e-mail: GEORG.SCHUMACHER@RUBA.RZ.RUHR-UNI-BOCHUM.DE\\

\noindent Pit-Mann Wong\\
Dept. of Mathematics, University of Notre Dame\\
Notre Dame, Indiana 46556\\
USA\\
e-mail: WONG.2@ND.EDU

\end{document}